\documentclass[11pt]{article}

\newcommand{\Hnko}{{\mathcal H^{k+1}}}
\newcommand{\Snkop}{{\mathcal S_+^{k+1}}}
\newcommand{\Snkp}{{\mathcal S^k_+}}

\usepackage{rotating}
\usepackage{algorithm2e}
\usepackage{algorithmic}
\usepackage{lineno}
\usepackage{datetime}
\usepackage{mathrsfs} 
\usepackage{ifthen}

\usepackage{kbordermatrix}
\usepackage{latexsym}
\usepackage{amsmath,amssymb,amstext,amsthm} 
\usepackage{graphicx,color,epsfig}
\usepackage{graphicx}
\usepackage{float}
\usepackage{verbatim}
\usepackage[bookmarks,pagebackref,
    pdfpagelabels=true, 
    ]{hyperref}
\usepackage{seceqn}
\usepackage{makeidx}
\makeindex

\usepackage{blkarray}
\usepackage{wrapfig,epic}
\usepackage{amscd}
\usepackage{url}

\DeclareMathAlphabet{\mathpzc}{OT1}{pzc}{m}{it}  

\newcommand{\GIF}{\mbox{\texttt{GIF}}}
\newcommand{\RRI}{\mbox{\texttt{RRI}}}
\newcommand{\gif}{\mbox{\texttt{GIF}}}
\newcommand{\MMtx}{\mbox{\texttt{M}}}
\newcommand{\gen}{\mbox{\texttt{gen}}}

\newcommand{\N}{\mathbb{N}}

\newcommand{\Symbol}{\mathcal{S}}

\newcommand{\Pro}{\widehat{\textbf{D}}}

\newcommand{\ProKer}{\textbf{D}}

\newcommand{\ProjKer}{\boldsymbol{\pi}}

\newtheorem{theorem}{Theorem}[section]   
\newtheorem{define}[theorem]{Definition}

\def\R{\mathbb{R}}
\def\H{\mathbb{H}}

\def\eqref#1{{\normalfont(\ref{#1})}}

\definecolor{rblue}{rgb}{.255,.41,.884} 

\def\eqref#1{{\normalfont(\ref{#1})}}

\newtheorem{example}{Example}[section]

\newtheorem{lem}{Lemma}[section]
\newtheorem{cor}{Corollary}[section]

\newtheorem{remark}{Remark}[section]

\newtheorem{problem}{Problem}[section]

\newcommand{\textdef}[1]{\textit{#1}\index{#1}}
\DeclareMathOperator{\nul}{{null}}
\DeclareMathOperator{\range}{{range}}

\newcommand{\Ss}{{\mathcal S} }

\newcommand{\LL}{{\mathcal L} }
\newcommand{\RR}{{\mathcal R} }
\newcommand{\EE}{{\mathcal E} }

\newcommand{\GG}{{\mathcal G} }

\newcommand{\PP}{{\mathcal P} }

\newcommand{\Scal}{{\mathcal S}}

\newcommand{\set}[1]{\left\{ #1 \right\}}

\newcommand{\Rm}{{\R^m}}

\newcommand{\Skp}{{\mathcal S^{k}_+}\,}

\newcommand{\A}{{\mathcal A}}

\newcommand{\bbm}{\begin{bmatrix}}
\newcommand{\ebm}{\end{bmatrix}}
\newcommand{\bem}{\begin{pmatrix}}
\newcommand{\eem}{\end{pmatrix}}
\newcommand{\beq}{\begin{equation}}
\newcommand{\beqs}{\begin{equation*}}
\newcommand{\bet}{\begin{table}}
\newcommand{\eeq}{\end{equation}}
\newcommand{\eeqs}{\end{equation*}}
\newcommand{\beqr}{\begin{eqnarray}}

\DeclareMathOperator{\trace}{{trace}}

\DeclareMathOperator{\svec}{{s2vec}}
\DeclareMathOperator{\sHvec}{{sHvec}}

\DeclareMathOperator{\sMat}{{s2Mat}}
\DeclareMathOperator{\sHMat}{{sHMat}}

\DeclareMathOperator{\rank}{{rank}}
\DeclareMathOperator{\spanl}{{span}}

\newcommand{\nc}{\newcommand}
\nc{\arrow}{{\rm arrow\,}}
\nc{\Arrow}{{\rm Arrow\,}}
\nc{\BoDiag}{{\rm B^0Diag\,}}
\nc{\bodiag}{{\rm b^0diag\,}}

\nc{\Mm}{{\mathcal M}^{m} }
\nc{\Mmn}{{\mathcal M}^{mn} }
\nc{\Mnr}{{\mathcal M}_{nr} }
\nc{\Mnmr}{{\mathcal M}_{(n-1)r} }
\nc{\kwqqp}{Q{$^2$}P\,}
\nc{\kwqqps}{Q{$^2$}Ps}

\nc{\notinaho}{(X,S)\in \overline{AHO}(\A)}
\nc{\inaho}{(X,S)\in AHO(\A)}

\newcommand{\bea}{\begin{eqnarray}}%
\newcommand{\eea}{\end{eqnarray}}%
\newcommand{\beas}{\begin{eqnarray*}}%
\newcommand{\eeas}{\end{eqnarray*}}%
\newcommand{\Int}{{\rm int\,}}

{}

\newcommand{\Hnp}[1][]{\,\mathbb{H}_+^{\ifthenelse{\equal{#1}{}}{n}{#1}}}
\newcommand{\Hn}[1][]{\,\mathbb{H}^{\ifthenelse{\equal{#1}{}}{n}{#1}}}
\newcommand{\Dn}[1][]{\,\mathbb{D}^{\ifthenelse{\equal{#1}{}}{n}{#1}}}

\begin{document}
\bibliographystyle{plain}
\title{
Facial Reduction and SDP Methods for Systems of Polynomial Equations
}
             \author{
\href{}{Greg Reid}
        \thanks{
Dept. Appl. Math., University of Western Ontario, London, Ontario,
Canada}
\and
\href{}
{Fei Wang}
        \thanks{
Dept. Appl. Math., University of Western Ontario, London, Ontario,
Canada}
\and
\href{http://orion.math.uwaterloo.ca/~hwolkowi/}{Henry Wolkowicz}%
        \thanks{
Department of Combinatorics and Optimization,
          University of Waterloo, Waterloo, Ontario N2L 3G1, Canada.
Research supported in part by The Natural Sciences and Engineering
                Research Council of Canada (NSERC) and 
the U.S. Air Force Office of Scientific Research (AFOSR).
}
\and
\href{}
{Wenyuan Wu}
	\thanks{Chongqing Key Lab.\;of Automated Reasoning and Cognition,
 CIGIT {\em Email:} wuwenyuan@cigit.ac.cn.\;
  Partly supported by cstc2013jjys0002 and West Light Foundation
  of the Chinese Academy of Science.}
}
\date{\today}
          \maketitle



\begin{abstract}
The real radical ideal of a system of polynomials with finitely many complex roots
is generated by a system of real polynomials having only real roots and free of multiplicities.
It is a central object in computational real algebraic geometry and important as a preconditioner
for numerical solvers.
Lasserre and co-workers have shown that the real radical ideal of real polynomial systems 
with finitely many real solutions can be 
determined by a combination of 
semi-definite programming (SDP) and geometric involution techniques.
A conjectured extension of such methods to positive dimensional polynomial 
systems has been given recently by Ma, Wang and Zhi.  

We show that regularity in the form of the
Slater constraint qualification (strict feasibility) fails for the
resulting SDP feasibility problems. Facial reduction is then a
popular technique whereby SDP problems that fail strict feasibility
can be regularized by projecting onto a face of the convex cone of 
semi-definite problems.

In this paper we introduce a framework for combining facial reduction 
with such SDP methods
for analyzing $0$ and positive dimensional real ideals of real 
polynomial systems.
The SDP methods are implemented in MATLAB and our geometric involutive 
form is implemented
in Maple.  We use two approaches to find a feasible moment matrix.
We use an interior point method within the CVX package for MATLAB and 
also the Douglas-Rachford (DR) projection-reflection method.

Illustrative examples
show the advantages of the DR approach for some problems over standard
interior point methods. We also see the advantage of
facial reduction both in regularizing the problem and also in
reducing the dimension of the moment matrices.
Problems requiring more than one facial reduction are also presented.
\end{abstract}


\section{Introduction}

In breathrough work Lasserre and collaborators \cite{LasserreLaurentRostalski09,MR2830310} 
have shown that the real radical ideal of real polynomial systems 
with finitely many real solutions can be 
determined by a combination of SDP and geometric involution techniques.
The real radical ideal of a system of polynomials with finitely many complex roots
is generated by a system of real polynomials only having real roots and free of multiplicities.
It is a central object in computational real algebraic geometry and important as a preconditioner
for numerical solvers.
A conjectured extension of such methods to positive dimensional polynomial 
systems has been given recently by Ma, Wang and Zhi \cite{MWZ:2012,MA12}.


The above approaches use the \textdef{method of moments} and the \textdef{Semi-definite
Programming, SDP} formulation.  In this paper we see that the Slater constraint
qualification, strict feasibility, fails for the SDP formulation
resulting in an ill-posed feasibility problem. 
Our main contribution is to use \textdef{facial reduction} to
project the problem onto the \textdef{minimal face} to help regularize these computations.
Our approach provides tools for working with the ideals involved, 
and gathering data on the open problem above.

\subsection{SDP and Facial Reduction}

The SDP formulation of the moment problem is equivalent to finding
$X$ for the linear feasibility system
\begin{equation}
\label{eq:sdpfeas}
  \A X = b, \quad X \in \Skp,
\end{equation}
where $\Skp$ denotes the convex cone of 
\index{$\Skp$, semi-definite cone}
\index{semi-definite cone, $\Skp$}
$k\times k$ real symmetric positive semi-definite matrices, and
$\A : \Skp \rightarrow \R^m$ is a linear transformation.
The standard regularity assumption for \eqref{eq:sdpfeas}
is the \textdef{Slater constraint qualification} or strict feasibility
assumption:
\begin{equation}
\label{eq:Slater}
\text{there exists } \hat X \text{ with }
  \A \hat X = b, \quad \hat X \in \Int \Skp.
\end{equation}
We let $X\succeq 0, \succ 0$ denote $X \in \Skp, \in \Int \Skp$,
respectively.
It is well known that the Slater condition holds generically,
e.g.,~\cite{DurJaSt:12}. Surprisingly,
many SDP problems arising from particular applications,
and in particular our polynomial system applications,
are marginally infeasible, i.e.,~fail to satisfy strict feasibility.
This means that the feasible set lies in the boundary of the cone, 
and even the slightest perturbation can make the problem infeasible.
This creates difficulties with the optimality and duality conditions as
well as with numerical algorithms.
To help regularize such SDP problems so that \textdef{strong duality} holds, 
facial reduction was introduced 
in 1982 by Borwein and Wolkowicz \cite{bw1,bw3}.
However it was only much later that the power of facial reduction was
exhibited in many applications, e.g.,~\cite{KaReWoZh:94,WoZh:96,AlWo:99}.
Developing algorithmic implementations of facial reduction that work
for large classes of SDP problems and the connections with
perturbation and convergence analysis has recently been achieved in
e.g.,~\cite{kriswolk:09,ChDrWo:14,ChWosensit:14,DrusLiWolk:14}.

A polynomial system of equations
can be viewed as a linear (or coefficient matrix) 
function of its monomials \cite{LasserreLaurentRostalski09,MR2830310}.
This linear function yields part of
the system of linear constraints in the SDP formulation of 
polynomial systems. The convex cone for polynomials are semi-definite 
moment matrices encoding the real solutions of the polynomial equations
and certain generalized Macaulay structure possessed by
the  polynomial systems.
Remarkable advances have been recently made in this area 
\cite{LasserreLaurentRostalski09,MR2830310,BPT2012} which
is an intersection between optimization and algebraic geometry.
In this article we establish a framework for using facial 
reduction for such systems and then solving the systems using the
regularized smaller SDP. 

\subsection{Prolongation projection methods for involutive bases of polynomial systems}

We now look at the details in the 
semi-definite linear constraint $\A X=b$ for the polynomial systems.
Polynomial systems are remarkable, in that many of 
their constraints are \textit{hidden}.
For example consider the degree two system 
\[
x^2 -  x - 1 = 0, \quad  x y - y - 1 = 0.  
\]
A single \textdef{prolongation} of this system to degree $3$ is found by 
multiplying them by each of the variables $x$ and $y$:
\begin{equation}
\label{eq:SimpExPro}
\begin{array}{rcl}
  x( x^2 - x - 1 ) &=&  x^3 - x^2 - x  \\
     x(x y - y - 1) &=&  x^2 y - x y - x \\
  y( x^2 - x - 1 ) &=&  x^2 y - x y - y \\
     y(x y - y - 1) &= & x y^2 - y^2 - y.
\end{array}
\end{equation}
\textdef{Projecting} in our paper loosely means eliminating higher degree monomials in favour of lower degree ones.
In the prolonged system we can project the system from degree $3$ to degree $2$ by eliminating the highest degree term
$x^2 y$ that occurs in the second and third equations of (\ref{eq:SimpExPro}):
\begin{equation}
\label{eq:SimpExProj}
  \left\{
\begin{array}{c}
x^2 y - x y - x = 0 \\   x^2 y - x y - y = 0  
\end{array}
\right\}
\implies xy + x = xy + y.
\end{equation}
Consequently we obtain the new projected (hidden) constraint $x = y$.
This process of uncovering the hidden polynomial constraints by prolongation and projection
 is effected numerically through our geometric involutive form
algorithm which has been implemented in Maple \cite{SRWZ2010,ReidLinWittkopf:2001}.

We note that
familiar methods for linear systems of equations are 
\textdef{Gaussian elimination, GE}, 
for exact solutions and \textdef{singular value
decompositions, SVD}, for least squares solutions.
For polynomial systems, the corresponding method in the exact case uses
\textdef{Gr\"obner Bases} \cite{BasuPollackRoy06}; while in the approximate
case we use \textdef{geometric involutive bases} \cite{SRWZ2010}.

\subsection{Facial Reduction and SDP methods applied to real radical
ideals of polynomial systems}

A major motivation for our paper is the success of the
work of Lasserre et al \cite{LasserreLaurentRostalski09} which gives
a new symbolic-numeric approach for computing the real radical ideal of 
zero dimensional polynomial
systems using geometric involution and SDP techniques.
Zero dimensional real polynomial systems are systems
with real coefficients and finitely many complex and real roots.
 Another major motivation is the important work on this topic in 
\cite{MWZ:2012,MA12} which 
conjectures
an extension of \cite{LasserreLaurentRostalski09} to positive 
dimensional real radical ideals.
Such ideals have associated real solution components (manifolds) of 
dimension $\geq 1$. (See also the paper \cite{ReidWangWu:14} for 
examples and many references.)

The \textdef{real radical ideal, \RRI}, of our system $P$ is the set of all
\index{\RRI, real radical ideal}
polynomials with the same zero set as $P$.
To give the reader an informal introduction to 
{\RRI}s and their interpretation,
consider the simple case of \textdef{univariate polynomials}
with real coefficients, $n=1$.
In particular, a real univariate polynomial $p(x)$ can be 
factored in real factors $(x - a_j )$
and conjugate complex factors $(x - \alpha_\ell )$, 
$(x - \bar{\alpha}_\ell )$ so that
\begin{equation}
\label{eq:univar-fact}
p(x) = \Pi_j (x - a_j )^{d_j}  \Pi_k  (x - \alpha_k )^{r_k} (x -
\bar{\alpha}_k)^{r_k},
\end{equation}
where $d_j$ and $r_k$ are the multiplicities of the roots.
The \textdef{real polynomial ideal} generated by $p(x)$ is the set of 
polynomials of the form $ g(x) p(x) $ where 
$g(x)$ is any real polynomial. The {\RRI} of $p(x)$ is generated by the polynomial
\begin{equation}
\label{eq:RealRadPol}
q(x) =  \Pi_j (x - a_j ).
\end{equation}
In many applications we are only interested in real roots, and the {\RRI}
shown here discards all the complex roots.
Moreover it also discards multiplicities which is important in improving conditioning for polynomial solvers.
Many general polynomial system solvers, that are capable of determining all 
solutions explicitly or implicitly, compute all complex and real roots first. 
In particular a generic system of $n$ degree $d$ polynomials in $n$ variables
generically has $d^n$ roots and potentially very few roots.
Thus the development of methods that 
avoid the calculation of the complex roots and multiplicities
is important for efficiency of polynomial system solvers.

\subsection{Outline}

Since we use sophisticated results from diverse areas,
in Section \ref{sec:setup} we present basic ideas and objects through simple 
examples. We give a preliminary introduction to moment matrices and also 
give a preliminary simple illustration of the power of facial reduction in 
Section \ref{sec:univarclass}.

In Section \ref{sec:GeoInvBases} we give a condensed and more formal 
description of geometric involutive bases and related algorithms.
In Section \ref{s:MM} we discuss moment matrices and related algorithms.

In Section \ref{sec:BackgroundProj} we discuss the methods we used to solve our 
SDP feasibility problems.
Since the polynomial problems we consider fail strict feasibility,
we will use facial reduction to regularize them.
However standard  primal-dual interior point
semi-definite programming packages do not deliver the accuracy required to 
guarantee facial reduction.  This motivates us to use 
\emph{Douglas-Rachford (DR)} projection/reflection methods.
\index{Douglas-Rachford, DR}
\index{DR, Douglas-Rachford}
\index{alternating projection, MAP}
\index{MAP, alternating projection}

In Section \ref{sect:facialimpl} we will discuss our implementation of facial reduction.
In Section \ref{sect:numerics} we give numerical experiments. Our
concluding remarks are in Section \ref{sect:Conclusion}.

\section{Basic setup and illustrative examples}
\label{sec:setup}

This paper uses sophisticated methods from diverse areas. To help the reader,
we informally introduce the methods of the paper and illustrate them by 
simple examples.  This helps emphasize that the operations underlying our 
approach are reasonably straightforward.  

\subsection{Real polynomial systems}
\label{sec:RealPolySys}

For background and references to real algebraic geometry
and semi-definite programming see
e.g.,~\cite{BasuPollackRoy06,BPT2012,SaVaWo:97,MR2830310,AnjosLasserre:11}.

We consider a (finite) system of $\ell$ polynomials 
\index{system of $\ell$ polynomials, $P$}
\index{$P$, system of $\ell$ polynomials}
$P = \{ p_1, ... , p_\ell \}  \subset \R[x_1,\ldots,x_n]= 
\R[x]$, where $\R[x]$ is the set of all polynomials with 
real coefficients in the $n$ variables $x = \begin{pmatrix}x_1& x_2& 
\ldots & x_n \end{pmatrix}^T$. We let \textdef{$d=\deg(P)$}
denote the \textdef{degree of
the polynomial system}, i.e.,~the maximum of the degrees of the
polynomials $p_j$ in $P$. The solution set or variety of $P$ is
\index{variety of $P$, $V_\mathbb{K}$}
\index{$V_\mathbb{K}$, variety of $P$}
\begin{equation}
\label{eq:V(P)}
V_\mathbb{K}  (p_1,...,p_\ell) = \{ x \in \mathbb{K}^n :  p_j(x) = 0,\;
\forall 1\leq j \leq  \ell \}.
\end{equation}
This is the \textdef{real variety of $P$}  if $\mathbb{K} = \R$ 
and the \textdef{complex variety of $P$} if $\mathbb{K}= \mathbb{C}$.
The real ideal generated by $ P  =  \{   p_1,\ldots ,p_\ell  \}  
\subset \R[x]$  is:
\begin{equation}
\label{eq:<P>}
\left\langle P  \right\rangle_\R = \left\langle  p_1,\ldots,p_\ell  
\right\rangle_\R  = \{ f_1 p_1 + \ldots + f_\ell p_\ell  : 
f_j \in \R[x],  \forall 1\leq j \leq  \ell \}.
\end{equation}
\textdef{Monomials} are denoted by $x^\alpha := x_1^{\alpha_1} 
\cdots x_n^{\alpha_n}$, where $\alpha \in \N^n$, 
$\N$ is the set of nonnegative integers,
\index{$\N$, nonnegative integers}
and the \textdef{degree of $x^\alpha$} 
is $|\alpha| := \|\alpha\|_1=\alpha_1 + \cdots + \alpha_n$. 
It is clear that the
degree of each monomial $|\alpha|\leq d$, the degree of the polynomial.
Then for appropriate coefficients $a_{k,\alpha}$, and for each $k$,
\begin{equation}
\label{eq:sorted}
\begin{array}{l}
\text{we sort by total degree of } |\alpha|
\text{ in nondecreasing order}\\
\text{with components of } \alpha 
\text{ sorted in lexicographic order}.
\end{array}
\end{equation}
We can rewrite the system of $\ell$ polynomials, $P$, as
\index{system of $\ell$ polynomials, $P$}
\begin{equation}
\label{eq:P=Ax}
 P  =  \left\{ \sum_{ | \alpha|  \leq d} \;   a_{k, \alpha} \: x^\alpha :  k = 
           1, \ldots , \ell   \right\} .  
\end{equation}
Throughout this paper, we use graded reverse lexicographic order, 
which orders first by degree and then by
reverse lexicographic order. This order respects the Cartan class of variables,
which is important in our numerical determination geometric features of polynomial
systems such as those in Definition \ref{def:symbol-class-test}.

\begin{define}[Coefficient matrix of $P$, $C(P)$]
\label{defCoeffMtx}
\index{coefficient matrix of $P$, $C(P)$}
\index{$C(P)$, coefficient matrix of $P$}
Let $\textbf{x}^{(\leq d )}$ be the column vector of monomials $x^\alpha$
with $0 \leq  | \alpha | \leq d$ sorted as in \eqref{eq:sorted}.
Suppose that the coefficients $a_{k,\alpha}$ in \eqref{eq:P=Ax}
are similarly sorted.
Then define the coefficient matrix of $P$ by $C(P) = (a_{k,\alpha})$. 
\end{define}
The following lemma follows immediately.
\begin{lem}
With $C(P),\textbf{x}^{(\leq d)}$ defined in Definition
\ref{defCoeffMtx}, we have 
\[
P = C(P) \textbf{x}^{(\leq d )},
\]
with $C(P) \in \R^{\ell \times N(n,d) }$ and 
\index{$N(n,d)$}
     $N(n,d) := 
    \small{\left( \begin{array}{c}
    d + n \\
       d \\
 \end{array} \right) }$ 
is the number of monomials in $\textbf{x}^{(\leq d )}$.
\end{lem}

The well-known presentation of polynomial systems as linear functions of their monomials and the related coefficient matrix and its kernel and rowspace has been exploited in \cite{Stetter:2004,Mourrain:1996,Mourrain:1999,MollerSauer:2000} 
and in the historical work by Macaulay \cite{Macaulay:1916}.

\begin{example}
\label{ex:Deg8Univar}
Consider the system of two univariate polynomials
\begin{equation}
\label{eq:DeguUnivar}
P = \{ x^8  -  x^4 - 2 ,   x^8  - 3 x^4  + 2 \}  \subset \R[x].
\end{equation}
Here the coefficient matrix is given by $C(P)$ in the equations
\begin{equation}
\label{CMtxSimpleEx}
C(P)  \mbox{\textbf{x}}^{(\leq 8 )}=    \left( \begin{array}{rcccccccr}
    -2 & 0 & 0 & 0 & -1 & 0 & 0 & 0 &  1 \\
     2 & 0 & 0 & 0 & -3 & 0 & 0 & 0 &  1 \\
\end{array} \right)
\small{ \left[ \begin{array}{c}
	1 \\
	x \\
	\vdots \\
	x^7 \\
	x^8  \\
\end{array} \right]
 }
=
\small{
 \left( \begin{array}{c}
        0   \\
	 0  \\
\end{array} \right)
}
\end{equation}
A familiar computation for many readers is to eliminate the polynomials 
using a Gr\"obner basis calculation:  
$ x^8  -  x^4 - 2 - (x^8  - 3 x^4  + 2) =
2 x^4 - 4$ or equivalently $x^4 - 2$. 
The original $8$ degree polynomials can be discarded since they are consequences
of $x^4 - 2$. 
In particular  $x^8  -  x^4 - 2 = x^4 (x^4 - 2 )  +   (x^4 - 2) = (x^4 + 1) ( x^4 - 2)$ so 
it lies in the ideal generated by $x^4 - 2$.
Similarly  $x^8  - 3 x^4  + 2$ lies in the ideal generated by $x^4 - 2$ and can be discarded.
All polynomials in the ideal generated by $P$ are polynomial 
multiples of the single polynomial
\begin{equation}
\label{GBsimpex}
 x^4 - 2.
\end{equation}
It is easy to see that every system of univariate polynomials is 
equivalent to a single univariate polynomial by applying such simple operations.
For systems of multivariate polynomials, such a minimal object is 
called a Gr\"obner basis.
Gr\"obner bases have been intensively studied \cite{Cox96} and usually consist of 
several polynomials.
We use the geometric involutive form algorithm
discussed in Section \ref{sec:GeoInvBases}
to obtain a numerically stable cousin of Gr\"obner bases.

\end{example}

\subsection{Moment matrices and polynomials }
\label{s:MM-polys}

Moment matrices combined with SDP provide a method to discard the complex roots in polynomial systems with finitely many roots, such as the two complex roots
of $x^4 - 2$ in Example \ref{ex:Deg8Univar} above.
Here we focus on the construction of moment matrices. For theoretical
background the reader is directed to 
e.g.,~\cite{AnjosLasserre:11,LaurentRostalski:12}.

A moment matrix is an infinite real symmetric matrix $M = (M_{\alpha,
\beta})$ with indices corresponding to the indices of the monomials 
$\alpha, \beta \in\N^n$.
Here $\alpha$ is the index for rows and $\beta$ is the index for
columns. Without loss of generality, we assume that $M_{0,0} = 1$.

\begin{define}[Moment matrix]
\label{def:momentM}
Let $u=\left\{ u_\alpha : \alpha \in \N^n, |\alpha| \leq d \right\} \in
\R^{\scriptsize{N(n,d)}}$ be a vector of
indeterminates where the entries are indexed corresponding to
the exponent vectors of the
monomials in $n$ variables of degree at most $d$. The degree $d$ moment
matrix of $u$ is a 
$N(n,d) \times N(n,d)$ symmetric matrix with rows and
columns corresponding to monomials in $n$ variables of degree at most
$d$, and defined as
\[
M_d(u) = M(u)=\begin{bmatrix} u_{\alpha+\beta}
\end{bmatrix}_{|\alpha|,|\beta|\leq d}.
\]
\end{define}

Given a multivariate polynomial system $P \subset 
\R[x]$,
with  $d = \deg(P)$ and $M  \in \R^{N(n,d) \times N(n,d)}$ be the 
truncated real symmetric moment matrix.
The linear constraints imposed by $P$ are, see \eqref{eq:uuT} below,
\[
C(P) \: M  = 0,
\]
where $C(P)$ is the coefficient matrix function given in Definition \ref{defCoeffMtx}.

\begin{example}[Moment matrix for univariate example $x=(x_1)$]
\label{ex:MMUnivar}

The moment matrix in the univariate ($n=1$) case is the infinite matrix whose $(\alpha, \beta)$ entry
is $u_{\alpha+\beta}$ and $\alpha, \beta \in \N$
given by:

\begin{equation}\label{MMtxInfinite}
M(u)=  \begin{bmatrix}
    u_0 & u_1 & u_2 & u_3 & u_4 & \cdots  \\
    u_1 & u_2 & u_3 & u_4 & u_5 & \cdots \\
    u_2 & u_3 & u_4 & u_5 & u_6 & \cdots \\
    u_3 & u_4 & u_5 & u_6 & u_7 & \cdots  \\
    u_4 & u_5 & u_6 & u_7 & u_8 & \cdots \\
    \vdots & \vdots & \vdots & \vdots & \vdots & \ddots \\
  \end{bmatrix}, \qquad u_0=1.
\end{equation}
Note that (\ref{MMtxInfinite}) is a \textdef{Hankel matrix}.
In Example \ref{ex:Deg8Univar} a degree $8$ input system was reduced to a degree $4$
output polynomial $P = \{ x^4 - 2 \} $.
 Let us associate $u_\alpha \leftrightarrow x^\alpha$.
Then we recover the polynomial equation using the coefficient matrix as
\[
     C(P)  \mbox{\textbf{u}}_{(\leq 4 )} =
 \left(
            \begin{array}{ccccc}
              -2 & 0 & 0 & 0 & 1 \\
            \end{array}
          \right)
          \left(
            \begin{array}{c}
              1 \\
              u_1 \\
              u_2 \\
             u_3 \\
             u_4 \\
            \end{array}
          \right)
          = 0. 
\]
This implies that in terms of the solution $x$:   
\begin{equation}
\label{eq:uuT}
     C(P)  \mbox{\textbf{x}}^{(\leq4)}  (\mbox{\textbf{x}}^{(\leq4)})^T =
 \left(
            \begin{array}{ccccc}
              -2 & 0 & 0 & 0 & 1 \\
            \end{array}
          \right)
          \left(
            \begin{array}{c}
               1 \\
              x  \\
              x^2 \\
              x^3\\
              x^4 \\
            \end{array}
          \right)
          \left(
            \begin{array}{c}
              1\\
              x \\
              x^2  \\
             x^3 \\
             x^4\\
            \end{array}
          \right)^T
          = 0. 
\end{equation}
In the SDP-moment matrix approach we impose $u_0 = 1$.
We note that the association $u_\alpha \leftrightarrow x^\alpha$ extends to the formal correspondence 
$x^\alpha x^\beta  \leftrightarrow u_{\alpha+\beta}$. This allows for the
construction of the truncated moment matrix 
to degree $d =4$ of the polynomial system as:
\begin{equation}\label{MMtx4}
M(u) = 
   \left(
  \begin{array}{ccccc}
      1   & u_1 & u_2 & u_3 & u_4  \\
    u_1 & u_2 & u_3 & u_4 & u_5  \\
    u_2 & u_3 & u_4 & u_5 & u_6  \\
    u_3 & u_4 & u_5 & u_6 & u_7  \\
    u_4 & u_5 & u_6 & u_7 & u_8  \\
  \end{array}
\right).
\end{equation}
Appending the linear constraints, we get
\begin{equation}
\label{eq:cons-univ}
C(P)  M(u) = 0 \: .
\end{equation}
The linear constraints (\ref{eq:cons-univ}) are:
\begin{equation}
\label{eq:cons-univ-explicit}
\{ u_4 - 2 = 0,
u_5 - 2 u_1 = 0 ,
u_6 - 2 u_2 = 0,
u_7 - 2 u_3 = 0,
u_8 -  2 u_4  = 0 \}
\end{equation}
which via the correspondence $u_\alpha \longleftrightarrow x^\alpha$ is 
equivalent to
$ \{ x^4 - 2, x^5 - 2x,  x^6 - 2x^2, x^7 - 2x^3, x^8 - 2 x^4 \}$.
The equivalent SDP problem here is to find a maximal rank generic point
$u = ( u_\alpha )$ where $| \alpha | \leq 2 d $ in the moment matrix with
\begin{equation}
\label{eq:SDPex}
M(u) \succeq 0 , \hspace{1.0cm} C(P) M(u) = 0 \: .
\end{equation}
By imposing these simple linear constraints we get an explicit simplified
moment matrix problem in only three variables:
\begin{equation}\label{MMtx4c}
M(u) =
\left[
  \begin{array}{ccccc}
     1   & u_1 & u_2 & u_3 & 2    \\
    u_1 & u_2 & u_3 & 2    & 2u_1 \\
    u_2 & u_3 & 2    & 2u_1 &2u_2 \\
    u_3 & 2    & 2u_1 & 2u_2 & 2u_3 \\
    2    & 2u_1 & 2u_2 & 2u_3 & 4   \\
  \end{array}
\right] \succeq 0.
\end{equation}
We note that the substitution of the linear constraints to simplify the
problem and reduce the number of variables is equivalent to facial
reduction; see Section \ref{sect:facialimpl} below.
This moment matrix problem in \eqref{MMtx4c} is then sent to an SDP solver
to approximately find a vector $(u_1, u_2, u_3)$ 
if possible such that
$M$ is a positive semi-definite matrix with maximum rank.
This solver returns an approximation which can be recognized for illustrative convenience
as $(u_1, u_2, u_3) = (0, \sqrt{2}, 0), u_0=1,u_4=2$.
Its associated moment matrix and moment matrix kernel are:
\[
\begin{array}{l}
	M =  
\begin{bmatrix} 1&0&\sqrt {2}&0&2\\ 
0& \sqrt {2}&0&2&0\\ 
\sqrt {2}&0&2&0&2\,\sqrt {2} \\ 
 0&2&0&2\,\sqrt {2}&0\\ 
2&0&2\, \sqrt {2}&0&4
\end{bmatrix},\\
\ker M =  \spanl_\R\left\{   
\begin{pmatrix}
-2\\ 0\\ 0\\ 0 \\ 1\end{pmatrix},
\begin{pmatrix} -\sqrt {2}
\\ 0\\ 1\\ 0
\\ 0\end{pmatrix}, \begin{pmatrix} 0\\ -\sqrt {2}
\\ 0\\ 1\\ 0
\end{pmatrix} \right\}.
\end{array}
\]
The kernel yields the generating set of three polynomials
\begin{equation} 
\label{sdpker}
\begin{array}{rcl}
\Scal &=&
 \{ -2 + x^4,-\sqrt{2} + x^2,   -\sqrt{2}x + x^3 \}
\\&=&
 \{ (\sqrt 2 + x^2)(-\sqrt 2 + x^2),-\sqrt{2} + x^2,   x(-\sqrt{2} +
x^2) \}.
\end{array}
\end{equation}
The factorization in \eqref{sdpker} allows a trivial
Application of the geometric involutive form algorithm that yields a 
geometric involutive basis 
\begin{equation} \label{genset}
\{ -\sqrt{2} + x^2  \}.
\end{equation}
The first and third polynomials in \eqref{sdpker}
are a consequence of $-\sqrt{2} + x^2 $ by our
inclusion test, so are discarded, e.g.,~\cite{LaurentRostalski:12}. 
Thus  we have a basis of the {\RRI} in \eqref{genset}.
There are efficient eigenvalue methods that can exploit this geometric form to efficiently numerically compute the
roots as eigenvalues \cite{ReidZhi2009,ReidTangZhi:2003,Stetter:2004,Mourrain:1999}.
For such solving methods tailored to the real radical and its advantages
see \cite{LasserreLaurentRostalski09}.
The degree $8$ system trivially has two real roots given by the polynomial in
\eqref{genset}, i.e.,~$\pm 2^{1/4}$.

\end{example}


\subsection{A class of univariate geometric polynomials}

\label{sec:univarclass}

In this section we experimentally explore the behavior of our facial reduction approach
(Facial Douglas-Rachford, or abbreviated as FDR) compared to a standard
SDP solver (Yalmip SDP, abbreviated as YSDP) which does not use facial
reduction.  In particular we consider the class of univariate geometric polynomials which are
the partial sums to odd degree $d$ of the geometric series:
\[
p_d (x) = 1 + x + x^2 + \cdots + x^{d-1} + x^d
\]
where $d = 1, 3, 5, \ldots$.  Then for odd degree $d$ we have
\[
p_d(x) = (x + 1)(1 + x^2 +  \cdots  + x^{d-3} + x^{d-1} )
\]
where the even degree factor $1 + x^2 +  \cdots  + x^{d-3} + x^{d-1} $ 
has only complex roots.
The $d$ roots are $x = \exp  \left( \frac{2 j \pi i }{d + 1} \right)$,
$j = 1, \cdots ,d$, and the non-real
roots appear in complex conjugate pairs.
Consequently a generator for the {\RRI} is $x+1$.\footnote{We denote the
generator of the {\RRI}  by
$\sqrt[\R]{ \langle p_d (x) \rangle_\R } = \langle x + 1 \rangle_\R$.}

We solved this class of problems for odd degrees $d$ using both
the FDR\footnote{The Facial reduction Douglas Rachford method is
presented in Section \ref{sect:FDR} below.} method with MATLAB R2013b 
and the YSDP (Yalmip SDP, R20140605) method.
We used a laptop (Windows 8.1, Intel Core(TM) i7-4600U CPU @2.10GHz 2.70 GHz, 8GB RAM, 64-bit OS, x64-based processor).

\begin{figure}[h!]
\label{Fig:YSDP-FDR-Univ1}
\includegraphics[width=1.0 \textwidth]{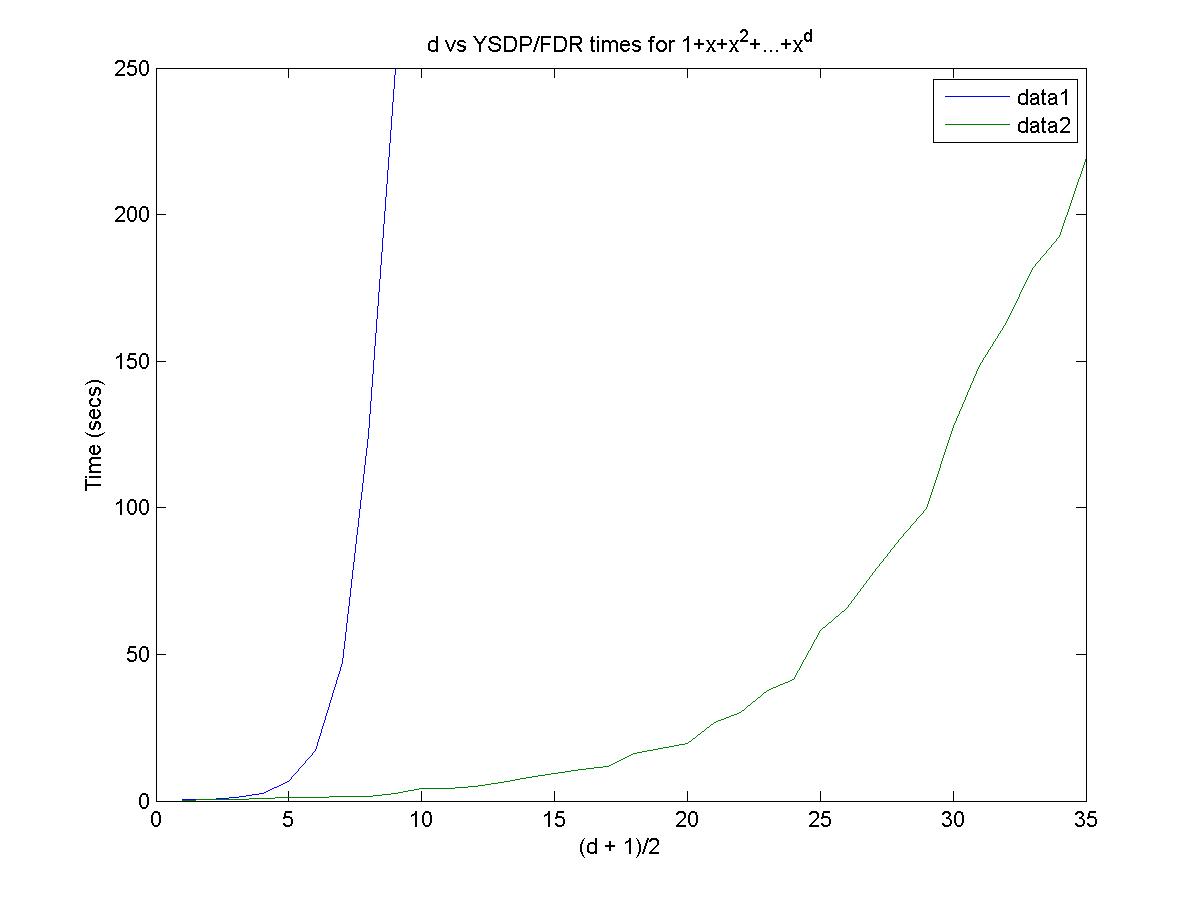}
\caption{Times (cpu secs) for the FDR method versus the YSDP methods 
applied to $p_d(x) = 1 + x + \ldots + x^d$ for odd degrees $1 \leq d \leq 69$.
The blue curve (data1 on the left) shows YSDP times and the green curve 
(data2 on the right) shows the significantly better FDR times.
}
\end{figure}

The running times (in cpu secs) for both methods are given in Figure 
\ref{Fig:YSDP-FDR-Univ1};
the range of values for the FDR method is clearly better. 

\section{Geometric involutive bases}
\label{sec:GeoInvBases}

In this section we introduce the basic objects for geometric involutive bases.
For details and examples see \cite{ReidWangWu:14,BLRSZ:2004}.

Involutivity originates in the geometry of differential equations.  See Kuranishi \cite{Kuranishi:1957} for a famous proof of termination of Cartan's prolongation algorithm
for nonlinear partial differential equations.
A by-product of these methods has been their implementation for linear
homogeneous partial differential equations with constant coefficients, and consequently for
polynomial algebraic systems.  See \cite{GerdtBlinkov:1998} for applications and symbolic algorithms for
polynomial systems.
The symbolic-numeric version of a geometric involutive form was first described and implemented in
Wittkopf and Reid \cite{ReidWittkopf:2001}.  It was applied to approximate symmetries of differential 
equations in \cite{BLRSZ:2004} and to polynomial solving in \cite{ReidZhi2009,ReidTangZhi:2003,SRWZ2010}.
See \cite{WuZhi12} where it is applied to the deflation of multiplicities in
multivariate polynomial solving.

\begin{define}
Let $P$ be  (as usual) a finite subset of $ \R[x]$ of degree $d$.
The $k$-th prolongation of system $P$ is
$\Pro^k  (P) =  \{  x^\alpha p :  0 \leq \deg (x^\alpha p ) \leq d+k,
\alpha \in \N^n, p \in P \}$.
\end{define}
For example $\Pro^k  (P)$ for $P = \{ x^2 - x - 1 , x y - y - 1 \}$
consists of $P$ together with the $4$ polynomials in (\ref{eq:SimpExPro}).  

\begin{define}
Given a subspace $V$ of $J^d := \R^{N(n,d)}$ and $\ell \leq d$,
define $\ProjKer^\ell (V)$ as the vectors of $V$ with
the components of degree $\geq d - \ell$ discarded.
Given $P \subset \R[x]$ of degree $d$
define $\ProjKer^\ell (P) := \ProjKer^\ell  \ker C(P)$.
The $k$-th prolongation of the kernel is
$\ProKer^k  (P) :=  \ker C(\Pro^k P)$.
\end{define}

See for example \cite{SRWZ2010} and the published references in 
\cite{ReidWangWu:14} for the stable numerical implementations of 
this paper's operations using SVD methods.
In Remark 3.5 of \cite{ReidWangWu:14} we discuss how prolongation and
projection can equivalently be computed in the kernel or rowspace, and
how polynomial
generators can always be extracted.  Underlying this is a 1 to 1 correspondence between the relevant vector spaces (not elements).


\begin{define}[\textbf{Symbol, class and Cartan involution test}]
\label{def:symbol-class-test}
Suppose $P \subset \R[x]$ of degree $d$.
The symbol matrix $\Symbol (P)$ of $P$ is the submatrix
of $C(P)$ corresponding to its degree $d$ monomials.
Then the class of a monomial $x^\alpha$ is the least $j$ such that $\alpha_j \not = 0$.
\end{define}


Suppose that the columns of $\Symbol (P)$ are sorted 
in descending order by class and that
it is reduced to Gauss echelon form.  For $k = 1, 2, ... , n$ define the quantities $\beta^{(k)}_d$ as the number of pivots in this reduced matrix of class $k$.  In a generic system of coordinates
the symbol is involutive if


\begin{equation}\label{CartanTest}
\sum_{k=1}^{k=n} k \beta_d^{(k)} = \rank \: \Symbol ( \Pro P ) 
\end{equation}


Suppose $Q \subset \R[x]$ has degree $d' $ and 
 a basis for $ \ker C(Q)$ is given by the rows of the matrix $B$.
To extract the $\beta_q^{(k)}$ in (\ref{CartanTest}) at projected degree
$d \leq d'$ we first numerically project $\ker C(Q)$ onto the subspace $J^d$
by deleting the coordinates in $B$ of degree $> d$ to give a spanning set $\tilde{B}$ 
 for  $\ProjKer^{d' - d} Q$.
Then delete the columns in $\widetilde{B}$
corresponding to variables of degree $< d$ to obtain a matrix $A_d$ corresponding to the
orthogonal complement of the degree $d$ symbol.
Let $A_d^{(k)}$ be the submatrix of $\widetilde{B}$ with columns corresponding to
variables of class $\leq k$.
In generic coordinates  for $ k = 1 \ldots  n $:
\begin{equation}  
\label{beta[k]}
\beta_d^{(k)}=   \small{ \left( \begin{array}{c}
							    n + d - k - 1 \\
							          d - 1 \\
 					\end{array} \right)  } -  \left( \rank \:  A_{d}^{(k-1)} - \rank \: A_d^{(k)} \right)  .
\nonumber
\end{equation}
Then the SVD can approximate the ranks in this equation for carrying out the Cartan Test 
  (\ref{CartanTest}).


\begin{define}
[Involutive System] A system of polynomials $P \subset \R[x]$ is involutive if
 $ \dim \:   \ProjKer \ProKer P = \dim \: P  $
and the symbol of $P$ is involutive.
\end{define}


\begin{define}
\label{def:projinvolsys}
Let $P \in \R[x]$ with $d = \deg P$ and
 $k$, $\ell$ be integers with $k \geq 0$ and $0 \leq \ell \leq k +d$.
 Then $\ProjKer^\ell \ProKer^k P$ is projectively involutive if
$ {\dim} \; \ProjKer^\ell \ProKer^k P  = {\dim}\;
\ProjKer^{\ell+1} \ProKer^{k+1} P$
and the { symbol} of $\ProjKer^\ell \ProKer^k P $ is involutive.
\end{define}
In \cite{BLRSZ:2004} we prove
that a system is projectively involutive if and only if it is involutive.
In the following algorithm we seek the smallest $k$ such that there exists an $\ell$
  with $ \ProjKer^\ell \ProKer^k P$ approximately involutive, and generates the same ideal as the input system. 
We choose the system corresponding to the largest such $\ell  \leq  k$ if there are several such values for the
given $k$.
\index{\GIF, Geometric involutive form}
\index{Geometric involutive form, \GIF}
\begin{algorithm}[h!]
\caption{{\GIF}: Geometric involutive form}
\label{alg:ProjInvBasis}
	\SetKwData{Input}{Input}
	\SetKwData{Output}{Output}
	\Input{
			   $ Q \subset \R[x_1,\ldots,x_n]$; tolerance $\epsilon$.
}\;
			Set $k := 0$, $ d := \deg(Q)$ and $P := \ker
C(Q)$\;
\While{$I\neq \emptyset$
}{
Compute $\ProKer^k (P)$;
initialize set of involutive systems $I:=  \{ \}$	\;
                	\For{$\ell$ from  $0$ to  $(d+k)$}{
 Compute $R := \ProjKer^\ell \ProKer^k (P)$\;
\If{$R$ involutive}{$I := I \cup \{ R \} $}
}
Remove systems $\bar{R}$ from $I$:  $\ProKer^{d + k-\bar{d}} 
  \bar{R} \not \subseteq \ProKer^k (P)$\;
					$k := k + 1$
}
\Output{  Return the polynomial generators  of the
{$\gif$}
 $(\bar{R})$ in $I$ of lowest degree $\bar{d} = \deg \bar{R}$.
}
\end{algorithm}

The degree of the geometric involutive basis in our method can be lower
than that given in \cite{MWZ:2012,MA12}
since Algorithm \ref{alg:ProjInvBasis} updates the generators with projections.
However in the absence of a proof of determination of the real radical the
larger moment matrices of \cite{MWZ:2012} can capture new members of the real radical in situations where our method has already terminated.

Additional discussion and examples are given in the long version of our work \cite{ReidWangWu:14}.

\section{Moment matrices \& algorithms}
\label{s:MM}

In this section we outline algorithms for combining geometric involutive form
and moment matrix methods; see Definition \ref{def:momentM}.
Recall that $M= M(u)=(M_{\alpha, \beta})$ denotes the moment matrix
indexed by $\alpha,\beta$ for rows and columns, respectively.
And, $d = \deg(P)$, $M  \in \R^{N(n,d) \times N(n,d)}$, and
the linear constraints imposed by our system of 
polynomials $P\subset \R[x]$ are given by the 
coefficient times moment matrix multiplication $ C(P)M =0$.
We let $\left\langle P \right\rangle_\R$ denote the
\textdef{associated polynomial ideal} and let
\[
\sqrt[\R]{ \left\langle P \right\rangle }_\R
              =
			 \{ f   \in \R[x] : f^{2m}+\sum_{j=1}^{s} 
q_j^{2} \in \left\langle P \right\rangle_\R, q_j
\in \R[x], m \in \N_+  \}.
\]
denote the 
\textdef{real radical ideal generated by polynomials $P$ over $\R$}.
A fundamental result \cite{BasuPollackRoy06} that is a consequence of the real nullstellensatz is
\[
\sqrt[\R]{ \left\langle P \right\rangle }_\R
              =  \{ f(x)   \in \R[x] : f(x) = 0,
\forall  x \in V_\R (P) \}.
\]

\begin{algorithm}[h!]
\caption{\texttt{GIF} -- \texttt{M} Method}
\label{alg:GIF-MMtx}
	\SetKwData{Input}{Input}
	\SetKwData{Output}{Output}
	\Input{
 $ P = \{p_1,...,p_k\} \subset \R[x_1,\ldots,x_n]$
}\;
Set $Q_0 := P,   \; \; \; j:= 0 $\;
\While{$r=d$
}
{
$  d :=  \dim \ker \GIF (Q_j)  $,  $ \; \; Q_{j+1} := \gen( \GIF (Q_j) )
$\;
Find $u^* = u( Q_{j+1}) \in \R^{N(n,2d)}$:  $M(u^*) \succeq  0,
C(Q_{j+1} ) M(u^*) = 0$\;
$  r := \rank( \MMtx (u^*) )$,  $ \; \;  Q_{j+2} := \gen(\ker \MMtx
(u^*) ) $\;
$j: = j + 2$
}
\Output{$Q_{j+1} \subset \R[x_1,\ldots,x_n]$;
  $ Q_{j+1} $ \textnormal{ is in geometric involutive form };
  $\sqrt[\R]{ \left\langle P
\right\rangle_\R  }  \; \;   \supseteq  \; \;  \left\langle Q_{j+1}
\right\rangle_\R  \; \;  \supseteq  \; \;  \left\langle P
\right\rangle_\R    $.
}
\end{algorithm}
Algorithm \ref{alg:GIF-MMtx} uses the following subroutines described as
Algorithms \ref{alg:MMtx-subroutine} and \ref{alg:gen-subroutine}.
\begin{algorithm}[h!]
\caption{\MMtx  \hskip3pt - Moment Matrix}
\label{alg:MMtx-subroutine}
	\SetKwData{Input}{Input}
	\SetKwData{Output}{Output}
	\Input{
$ Q \subset \R[x_1,\ldots,x_n]$.  \textnormal{Set} $ d := \deg(Q)$.
}\;
   Construct the moment matrix to degree $2d$.\;
   Use SDP methods to numerically solve for a generic point $u^* = u(Q)$ 
that maximizes the rank of the moment matrix subject to the constraints
$C(Q) \: M(u^* ) = 0$.\;
\Output{ \textnormal{ Return $\MMtx (u^*) \succeq 0$ the moment matrix evaluated at this generic point.}
}
\end{algorithm}

\begin{algorithm}[h!]
\caption{\gen}
\label{alg:gen-subroutine}
	\SetKwData{Input}{Input}
	\SetKwData{Output}{Output}
	\Input{
 $\GIF (Q)$ or $\ker \MMtx (u^*) $ where $u^* = u(Q)$. 
}\;
 \Output{\textnormal{Polynomial generators corresponding to $\GIF (Q)$ or $\ker \MMtx (u^*)$}
}
\end{algorithm}

\noindent

\begin{remark}[\textbf{Rank-Dim-Involutive Stopping Criterion}]
A natural termination criterion used in Algorithm \ref{alg:GIF-MMtx} is that the generators stabilize at some iteration and the system
is involutive:
\begin{equation}
\label{ConjMMtx}
\gen (\GIF (Q)) = \gen (\ker \MMtx (u^*)) \; \mbox{and} \;  Q \; \mbox{involutive where  } u^* = u(Q)
\end{equation}
By \cite{LasserreLaurentRostalski09} $\langle \gen( \ker \MMtx(Q_{j+1})) \rangle $ is a sequence of ideals containing 
$\sqrt[\R]{ \left\langle P \right\rangle }$ . 
We get an ascending chain of ideals in a Noetherian ring
$\R[x_1, ... , x_n ]$.  Hence, together with the finiteness of the Cartan-Kuranishi
geometric involutive form algorithm,  Algorithm \ref{alg:GIF-MMtx} terminates.
\end{remark}

\section{Mathematical background for the projection methods}
\label{sec:BackgroundProj}

In this section we describe the background for the projection methods
for finding feasible solutions for the moment problems.
An important part of these methods is building an efficient
matrix representation for the linear constraints
on the moment matrices resulting from the polynomial systems.


\subsection{Linear constraints for multivariate polynomial moment matrices}

Recall that we introduced moment matrices informally by a simple example
in Section \ref{s:MM-polys}; see also Definition \ref{def:momentM}.
Let $u_{\alpha} := u_{\alpha_1,...,\alpha_n}$ where $\alpha \in \N^n$
and  the degree of $u_{\alpha}$ is $| \alpha| = \alpha_1 + \ldots + \alpha_n$. 
Let $\langle {\alpha}_{(\leq d )} \rangle$ be an array of the subscripts 
$\alpha$ of $\langle u_\alpha \rangle$ with $0 \leq | \alpha| \leq d$
and sorted  as in \eqref{eq:sorted}.

Consider a truncated moment matrix $M(u) = (u_{\alpha+\beta})_{\alpha,
\beta \in \R^{N(d,n)}}$. 
The generalized truncated moment matrix can be represented as follows,
where $\langle \cdot \rangle$ yields the addition of the 
subscripts for the $f_j$:
\[
	M(u)  =\begin{bmatrix}
 \langle f_0(u),f_0(u) \rangle & \langle f_0(u),f_1(u) \rangle & \langle f_0(u),f_2(u) \rangle & \ldots & \langle f_0(u),f_l(u) \rangle \cr
\langle f_1(u),f_0(u) \rangle & \langle f_1(u),f_1(u) \rangle & \langle f_1(u),f_2(u) \rangle & \ldots & \langle f_1(u),f_l(u) \rangle \cr
\langle f_2(u),f_0(u) \rangle & \langle f_2(u),f_1(u) \rangle & \langle f_2(u),f_2(u) \rangle & \ldots & \langle f_2(u),f_l(u) \rangle \cr
\vdots & \vdots & \vdots & \ddots & \vdots \cr
\langle f_l(u),f_0(u) \rangle & \langle f_l(u),f_1(u) \rangle & \langle f_l(u),f_2(u) \rangle & \ldots & \langle f_l(u),f_l(u) \rangle \cr
\end{bmatrix} .
\]
Here, $\langle f_0,f_1,...,f_l \rangle$ corresponds to the array 
$\langle u_\alpha \rangle$ with $0 \leq | \alpha| \leq d$ sorted 
as in \eqref{eq:sorted}.
We denote the \textdef{i-th} element in 
$\langle u_\alpha \rangle$ by $u^i_{\alpha}$. Then $f_i(u)$ is $u^i_{\alpha}$.

In the univariate case the moment matrices have Hankel structure as shown in (\ref{MMtx4}).
In Table \ref{eq:BiMtx} we display a truncated bivariate moment matrix
partitioned into block submatrices having the same degree.
\begin{table}[h!]
\[
\small{
 M(u)=\left[ \begin {array}{c|c|c|c} 
   \begin{array}{c}u_{{00}}\end{array}&  
        \begin{array}{cc}u_{{10}}&u_{{01}}\end{array}  &
        \begin{array}{ccc} u_{{20}}&u_{{11}}&u_{{02}}\end{array}&
      \begin{array}{cccc} u_{{30}}&u_{{21}}&u_{{12}}&u_{{03}}\end{array} \\ 
\hline 
        \begin{array}{c}u_{{10}}\\u_{{01}}\end{array}  &
\begin{array}{cc} u_{{20}}&u_{{11}} \\ u_{{11}}&u_{{02}}\end{array}
& 
\begin{array}{ccc} u_{{30}}&u_{{21}} &u_{{12}} \\ u_{{21}}&u_{{12}} &u_{{03}} \end{array}
&
      \begin{array}{cccc} u_{{40}}&u_{{31}}&u_{{22}}&u_{{13}}\\
      u_{{31}}&u_{{22}}&u_{{13}}&u_{{04}} \end{array} \\ 
\hline 
        \begin{array}{c}u_{{20}}\\u_{{11}} \\u_{{02}} \end{array}  &
\begin{array}{cc} u_{{30}}&u_{{21}} \\u_{{21}}&u_{{12}} \\u_{{12}}&u_{{03}} \end{array}
& 
\begin{array}{ccc} u_{{40}}&u_{{31}} &u_{{22}} \\u_{{31}}&u_{{22}}&u_{{13}}  \\u_{{22}}&u_{{13}}&u_{{04}} \end{array}
&
      \begin{array}{cccc} u_{{50}}&u_{{41}}&u_{{32}}&u_{{23}}
      \\u_{{41}}&u_{{32}}&u_{{23}}&u_{{14}}  \\u_{{32}}&u_{{23}}&u_{{14}}&u_{{05}} \end{array} 
\\ \hline
        \begin{array}{c}u_{{30}}\\u_{{21}} \\u_{{12}}  \\u_{{03}} \end{array}  &
\begin{array}{cc} u_{{40}}&u_{{31}} \\u_{{31}}&u_{{22}} \\u_{{22}}&u_{{13}}  \\u_{{13}}&u_{{04}} \end{array}
& 
\begin{array}{ccc} u_{{50}}&u_{{41}} &u_{{32}} \\u_{{41}}&u_{{32}}&u_{{23}}  \\u_{{32}}&u_{{23}}&u_{{14}}  \\u_{{23}}&u_{{14}}&u_{{05}} \end{array}
&
      \begin{array}{cccc} u_{{60}}&u_{{51}}&u_{{42}}&u_{{33}}
      \\u_{{51}}&u_{{42}}&u_{{33}}&u_{{24}}  \\u_{{42}}&u_{{33}}&u_{{24}}&u_{{15}}  \\u_{{33}}&u_{{24}}&u_{{15}}&u_{{06}} \end{array} 
\end{array}
 \right] 
}
\]
\caption{
A truncated bivariate moment matrix
partitioned into block submatrices having the same degree.
}
\label{eq:BiMtx}
\end{table}
Notice that the matrix in Table \ref{eq:BiMtx} is not Hankel. 
However each of its block matrices is rectangular Hankel;
though even this feature is lost for multivariate moment matrices in more than two variables.

As mentioned above, without loss of generality we assume that
$u_{00} = 1$. As an abbreviation, we may denote $M=M(u) = M_d(u)$. 

Besides being a symmetric matrix, the moment matrix also has other 
linear constraints among its entries. 
One can easily see these constraints in the truncated univariate matrix  
(\ref{MMtx4}) and bivariate matrix in Table \ref{eq:BiMtx}.
An important requirement of our projection methods is to maintain 
these constraints.
For example, in the bivariate case above, the matrix elements
$M(u)_{14} =M(u)_{22} = u_{20}$ are equal.

We now outline a simple algorithm to find a non-redundant 
matrix representation of these constraints.
To list these constraints we start from the first row 
and traverse the matrix from left to right across the rows
and then traverse the rows from top to bottom.
Note also that we only need examine entries above the main diagonal 
since the matrix is symmetric.

For (\ref{MMtx4}) the first linear constraint traversing from the 
first row downwards is $M(u)_{14} = M(u)_{22}$.   
We denote $e_i$ as the \textdef{$i$-th unit vector} and 
$E_{ij} = \frac{1}{2}(e_i^T e_j + e_j^T e_i)$.
To impose this constraint, we construct matrix $A_t = E_{22} - E_{14}$,
where $t$ represents the index of the linear constraints 
and $t = 2$ in this case. The constraint is then given by
\[
\langle A_t,M\rangle= \trace( (E_{22} - E_{14})M)=0.
\]
Since we always assume $M(u)_{1,1} = 1$,  we need to set $A_1 = E_{11}$. 
Here $A_t$ is called the \textdef{matrix representative} of the
\textdef{t-th} linear constraint.
The collection of all such matrix representatives for a given moment matrix 
is called the \textdef{matrix representation} of the moment matrix structure.

Algorithm \ref{alg:linearconstraints} below determines all 
the (non-redundant) matrix representatives of the
linear constraints defining the
matrix representation of the multivariate moment matrix structure.










\begin{algorithm}[h!]
\caption{
Matrix representation of moment matrix structure
\label{alg:linearconstraints}
}
	\SetKwData{Input}{Input}
	\SetKwData{Output}{Output}
\Input{$d$, $n$}\;
\textnormal{Initialize array $T = \langle \alpha_{(\leq d )} \rangle $ and $T(i)$ is the \textdef{i-th} element of $T$}.  \\
\textnormal{Initialize n array $S = \langle s \rangle$ with the same length as $\langle \alpha_{(\leq d )} \rangle$ and $S(i) =[(1,i);\alpha_{(\leq d )}(i)]$ where $S(i)$, $\alpha_{(\leq d )}(i)$ is the  
\textdef{i-th} element of $S$, $\langle \alpha_{(\leq d )} \rangle$}.

 \textnormal{Let $m$ be the length of $T$, $t =2$ and $A_1 = E_{11}$}. 

	\For{$i$ from $2$ to $m$,}{
 		\For{ $j$ from $i$ to $m$,}{
			\eIf{
     there exists an $s =[(g,h);\alpha] \in S$ such that $T(i) + T(j) = \alpha $
                     }
      {$A_t = E_{ij}-E_{gh}$, $t = t +1$}
                {
 Adjoin a new element $s = [(i,j); \alpha]$ to $S$ where $\alpha = T(i) + T(j)$
			}
			}
		}
\Output{
\textnormal{Return an array of matrix representatives $\{A_t\}$ 
where $t \in \EE$, $\EE = \{1,2,\dots,\eta \}$ and $\eta$ is the total number of the linear constraints.}
}\;
\end{algorithm}
There are no redundant relations produced by this algorithm so we can avoid an overdetermined system.

In what follows for applications to multivariate polynomial systems of degree $d$ in $n$ variables we have 
\begin{equation}
\label{eq:def-k}
k :=   
N(n, d) =
\small{\left( \begin{array}{c}
  d + n \\
    d \\
 \end{array} \right) }
\end{equation}
Our main problem is the following.
\begin{problem}[Main Problem]
\label{prob:mainM}
Let $B$ be a given $(k+1) \times m$ matrix of full column rank. Find $u
\in \R^{2k+1}$ so that
\[
	B^TM(u)=0, \quad \trace E_{11}M(u) = 1, \qquad  M(u) \succeq 0.
\]
\end{problem}

We denote \textdef{$\Hnko$, space of generalized Hankel matrices}.
That is these matrices have the multivariate structure whose matrix representation 
is computed by Algorithm \ref{alg:linearconstraints}.
It is well known that the special case of Hankel matrices are notoriously ill-conditioned.
This means that the cone $\Snkop \cap \Hnko$ is {\em thin}, i.e.,~it is close
to the boundary of $\Snkop$, e.g.,~\cite{MR929571,MR1771780,MR1159043}.  
Therefore, solving Problem \ref{prob:mainM}
using semi-definite programming techniques results in
numerical difficulties.

\subsection{Methods of alternating projection and Douglas-Rachford
projection-reflection}
\label{sect:apmdr}

To apply the methods of \textdef{alternating projection, MAP} or 
\textdef{Douglas-Rachford reflection-projection}, we want to express the 
main Problem \ref{prob:mainM}
as an equivalent problem with moment matrix $M=M(u)$:
\begin{equation}
\label{eq:mainA}
\A(M)=b, \quad B^TM=0, \quad M \in \Ss_+^{k+1}.
\end{equation}
Here the linear transformation $\A$ is obtained from Algorithm
\ref{alg:linearconstraints}.
The following Corollary \ref{cor:solvesys} 
provides the details of the system that we want
to solve. We first apply facial reduction and get a smaller system.
Recall from Algorithm \ref{alg:linearconstraints}, we get an array 
of representing matrix $A_t$\thinspace s where $t \in \EE$, 
$\EE = \{1,2,\dots,\eta \}$.
\begin{cor}
\label{cor:solvesys}
Let $V$ be $(k+1) \times (k+1-m)$ and satisfy $V^TV=I, V^TB=0$.
Let $\bar A_{t} \leftarrow V^TA_{t}V, \forall t\in  \EE$. 
Let $\bar \A: \Ss^{k+1-m} \rightarrow \R^{ \EE}$ be defined by 
\begin{equation}
\label{eq:Abar}
	\bar \A(\bar M):=	\begin{pmatrix}  
		\left( \trace \bar A_{t}\bar M\right)_{t} \cr
	\end{pmatrix}_{\forall t \in  \EE}
\end{equation}
Then the main Problem
\ref{prob:mainM} with $V\bar{M}V^T=M(u)$ is equivalent to
\eqref{eq:mainA}, i.e.,~to
\[
\bar \A(\bar{M})=e_1, \qquad \bar M \in \Ss^{k+1-m},
\]
and we get $M(u)=V \bar{M}V^T$.
\qed
\end{cor}

Let $L$ denote the matrix representation for $\bar \A$ in
the linear constraints in Corollary \ref{cor:solvesys}.
There are two projections we use to update the current point $p_c$.
First, we look at \textdef{$\PP_{\LL}$, the linear manifold projection}.
For the linear system 
$Lp = b=e_1$ where $L$ has full row rank, we solve the nearest point problem
$\min \left\{\frac 12\|p-p_c\|_2^2 : Lp=b \right\}$, i.e.,~we find the
projection onto the linear manifold for the linear constraints.
We use \textdef{$L^{\dag}$, the Moore-Penrose generalized inverse} of $L$. 
The residual and the update $p_+$ are then 
\begin{equation}
\label{eq:resupdate}
r_c = b-Lp_c; \qquad p_{+} = p_c + L^{\dag}r_c.
\end{equation}
Second, we project the updated symmetric matrix 
$P_+=\PP_{\LL}(P_c)=\sHMat(p_+)$ onto the 
semi-definite cone using the 
Eckart-Young Theorem \cite{EckartYoung39}, i.e.,~we diagonalize and zero
out the negative eigenvalues. Here $\sHMat=\sHvec^*=\sHvec^{-1}$ is both the
adjoint and the inverse mapping. We denote
\textdef{$\PP_{\Snkp}$, the positive semi-definite projection} and get the
new positive semi-definite approximation $\PP_{\Snkp}(P_+)$.

\subsubsection{Method of alternating projections}
\index{alternating projection, MAP}
\index{MAP, alternating projection}
The MAP method is particularly simple, see e.g.,~the recent book
\cite{MR2849884}.
We begin with an initial estimate, e.g.,~$P_c = \alpha I \in {{\mathcal M}^{mk} } $ 
for a large $\alpha > 0$. We then repeat the projection steps in Items
\ref{item:projL}, \ref{item:projP}, \ref{item:update} till a
sufficiently small desired tolerance is obtained in the norm of the residual.
\begin{enumerate}
\item
\label{item:projL}
Evaluate the residual $r_c = b-Lp_c$. Use the residual to
evaluate the linear projection and obtain the update
\[
P_L=\PP_{\LL}(P_c).
\]
\item
\label{item:projP}
Evaluate the positive semi-definite projection
using the Eckart-Young Theorem and update the current approximation
\[
P_{SDP}=\PP_{\Snkp}(P_L).
\]
\item
\label{item:update}
Update the cosine value in \eqref{eq:cosvalue}. Then update
$P_c=P_{SDP}$.
\end{enumerate}
The (linear) convergence rate is measured using cosines of angles from three
consecutive iterates
\begin{equation}
\label{eq:cosvalue}
\cos (\theta)=  \left(\frac{\trace \left((P_L-P_{c})^*(P_{SDP}-P_L)\right)}
 {\left\|P_L-P_{c}\right\|\left\|P_{SDP}-P_L)\right\|}\right).
\end{equation}

\subsubsection{Douglas-Rachford reflection method}
\label{sect:FDR}
\index{Douglas-Rachford, DR}
\index{DR, Douglas-Rachford}
Recall the projections defined above $\PP_{\LL}, \PP_{\Snkp}$.
We want to find, see \eqref{eq:mainA},
\[
   P\in \GG\cap\Snkop, \quad\text{where}\ 
      \GG:=\set{P\in\Snkop : \A(P)=b}.
\]
We now apply the Douglas-Rachford (DR) 
projection/reflection method \cite{MR0084194}. 
(See also e.g.,~\cite{ArtachoBoTa:13,MR3149115}.)

Using the QR algorithm applied to $B$ and $A$, 
we start with an initial estimate 
\begin{equation}
   \label{eq:initP}
P_0 \succeq 0 \text{ with } B'P_0=0 \text{ and } (1,1) \text{ component } =1.
\end{equation}
Define the \emph{reflections} $\RR_{\LL},\RR_{PSD}:\Snkop\to\Snkop$
using the corresponding projections, i.e.,~\index{reflections, 
$\RR_{\LL},\RR_{PSD}$}
\index{$\RR_{\LL},\RR_{PSD}$, reflections}
\[
   \RR_{\LL}(P) := 2\PP_{\LL}(P) - P,\quad
   \RR_{PSD} (P) := 2\PP_{PSD}(P) - P,\quad
   \forall\, P\in\H^{mk}.
\]
\begin{itemize}
\item
\underline{\bf Initialization:}
We set our current estimate $P_c=P_0$ to satisfy \eqref{eq:initP}.
We calculate the residual $Res_{\LL}=R-A*\sMat(P_c)$, set
$normres=\|Res_{\LL}\|$,
denote the reflected residual $Resrefl_{\LL}=Res_{\LL}$ and reflected
point $\RR_{PSD}=P_c$.
\item
\underline{\bf Iterate:}
We continue iterating from this point
while $normres>toler$, our desired tolerance.
\item
We use Resrefl to project the current reflected PSD point 
$\RR_{PSD}$ onto the linear manifold to get the projected
point $P_{\LL}=\RR_{PSD}+A^\dagger Resrefl$.
Then we reflect to get our second reflection point
$\RR_{\LL}=2*P_{\LL}-\RR_{PSD}$
\item
At this time we set our new/current estimate for 
convergence to be $P_c= P_{new}=(P_c+\RR_{\LL})/2$.
\item
We now project $P_c$ to get $P_{PSD}$. We check the residual here
for the stopping criteria $normres=\|Res_{\LL}\|=\|R-\A P_{PSD}\|$.
\item
We now calculate the first reflection point $\RR_{PSD}=2*P_{PSD}-P_c$
and update the reflected residual $Resrefl=R-A \svec(\RR_{PSD})$.
\end{itemize}

The Douglas-Rachford projection/reflection method is simply:
\begin{enumerate}
\item
   Start at an initial point $P_0\in\Snkop$ satisfying \eqref{eq:initP}
\item
   Iterate: $P_{j+1} = \frac{1}{2} (P_j + 
      \RR_{PSD}(\RR_{\LL}(P_j))$, for all $j=0,1,\ldots$.
\end{enumerate}


Also the basic theorem on the convergence of
the sequence ${ \Pi_G (X_k) }_k$ ,
\cite[Thm 3.3, Page 11]{MR3149115},
\href{http://carma.newcastle.edu.au/jon/cycDRinfeas.pdf}
{carma.newcastle.edu.au/jon/cycDRinfeas.pdf}.
so the residuals of the projections of the iterates on one of the sets
have to be used for the stopping criteria.
We use the residual after the projection onto the SDP cone since finding
the residual with respect to the linear manifold is inexpensive.

To check the linear convergence rates we use the cosine of the angles for the
vectors of successive iterates, i.e.,~for three successive iterates
$P_c,\RR_{PSD},\RR_{\LL}$, and
\[
\cos (\theta) =  \left|\frac{\trace
\left((\RR_{PSD}-\RR_{\LL})^*(\RR_{PSD}-P_c)\right)}
 {\left\|(\RR_{PSD}-\RR_{\LL})\right\|\left\|\RR_{PSD}-P_c\right\|}\right|.
\]

\section{Facial reduction implementation}
\label{sect:facialimpl}
Our moment problem is a feasibility problem of the form
\begin{equation}
	\label{eq:elemmomentprob}
	B^T M(u) = 0, \quad M(u)\succeq 0,
\end{equation}
where $B$ is a given matrix and $M(u)$ is a linear function of the
variables $u$.
Constraints on $M(u)$ are described in Section \ref{sect:apmdr}, where
the problem is changed to equality form and then
uses facial reduction to get the form
\begin{equation}
\label{eq:mainAA}
\bar \A(P)=\bar b, \qquad  P \succeq 0.
\end{equation}
This form  includes the first step of facial reduction using the matrix
$B$, see Corollary \ref{cor:solvesys} and \eqref{eq:Abar}.
Here $\bar \A(P) = (\trace \bar A_i P)\in \Rm$, for specific symmetric
matrices $\bar A_i$.

The projection methods behave poorly when Slater condition fails. We
therefore attempt to apply further steps of facial reduction and reduce system
\eqref{eq:mainAA}  until a strictly feasible point exists.
We use the following theorem of the alternative or characterization
of a strictly feasible point; see e.g.,~\cite{ChWosensit:14}.
\begin{equation}
	\label{eq:thmalt}
	\begin{array}{cc}
	\exists \hat P, \bar \A(\hat P)=\bar b, \hat P \succ 0
	\\  \iff
	\\ Z=\bar \A^*y \succeq 0, \bar b^Ty=0 \implies Z=0.
\end{array}
\end{equation}
Note that if a $Z\neq 0$ can be found satisfying the left part of the
bottom half of \eqref{eq:thmalt} and for the top half 
$\hat P \succeq 0, \bar(\hat P)=\bar b$, then
\[
	0= \bar b^Ty = \langle \bar A (\hat P),y\rangle = \langle \hat
	P, Z \rangle \implies \hat PZ=0 \implies \range \hat P \subseteq
	\nul Z.
\]
Therefore, if the full column rank matrix $W$ satisfies
$\range W = \null Z$, then we can facially reduce the problem
using the substitution $\hat P=W \bar P W^T$, i.e.,~we can restrict the
feasibility problem  in \eqref{eq:mainAA} to the face $W\cdot W^T$.

We can implement the test in \eqref{eq:thmalt} in several ways.
We suppose that $\bar A$ is the matrix representation of $\bar \A$,
i.e.,~we let $p=\svec(P)$ and then we have
\[
	\bar A p =(\bar \A \sMat)( \svec (P)) = \bar \A P, \quad  \bar \A^* y 
	= \sMat(\bar A^T y).
\]
One way would be to first evaluate the orthogonal matrix
$\begin{bmatrix} \frac 1{\|b\|} b &U \end{bmatrix}$ and find $v$ so that
\[
	\sMat (\bar A^T (Uv)) \succeq 0, \quad 
	\trace \bar \A^* (Uv) =(\bar \A(I)^TU) v=1.
\]
Alternatively, we solve
	\footnote{
This can be implemented in e.g.,~CVX using the \emph{norm} function 
or absolute value function for the objective, i.e.,~we minimize $|\bar
b^Ty|$ rather than using the squared term.
}
\[
	\begin{array}{rcl}
		p^*
		:= &	\min & \frac 12 (\bar b^T y)^2 \\
		     &\text{s.t.} 
		& \bar \A^* y\succeq 0\\
		&&    \trace \bar \A^* y = 1\\
	\end{array}
\]

\section{Numerical experiments}
\label{sect:numerics}

\subsection{Examples of Ma, Wang and Zhi \cite{MWZ:2012}}

Ma, Wang and Zhi \cite{MWZ:2012,MA12} present an approach using Pommaret Bases coupled with moment matrix completion
to approximate the real radical ideal of a polynomial variety. 
We applied our approach to \cite[Examples 4.1-4.6]{MWZ:2012}.
with the results shown in Table \ref{tab:FDR-GIF}. 
In each case we obtained a geometric involutive basis which can be independently verified
as a geometric involutive basis for the real radical.
In \cite{MWZ:2012} Pommaret bases are successfully obtained for the real radical for these examples.

Here are the $6$ systems of polynomials corresponding to the examples in
\cite{MWZ:2012}:
\begin{subequations}
\begin{align}
 &  \{  x_1^2 + x_1 x_2 - x_1 x_3 - x_1 - x_2 + x_3 ,  \;\;  x_1 x_2 + x_2^2 - x_2 x_3 - x_1 - x_2 + x_3 ,     \nonumber  \\
\label{eq:MWZ41}
 &   \hspace{5cm} x_1 x_3 + x_2 x_3 - x_3^2 - x_1 - x_2 + x_3     \}     \\
\label{eq:MWZ42}
 &  \{ x_1^2 - x_2 , \; \;   x_1 x_2 - x_3 \}   \\
\label{eq:MWZ43}
 &   \{ x_1^2 + x_2^2 + x_3^2 -2  , \; \;   x_1^2 + x_2^2 - x_3 \}   \\
\label{eq:MWZ44}
 &  \{ x_3^2 + x_2 x_3 - x_1^2 ,    \;\;            x_1 x_3 + x_1 x_2 - x_3,  \;\;   x_2 x_3 + x_2^2 + x_1^2 -x_1  \}       \\
\label{eq:MWZ45}
 &  \{ (x_1 - x_2)(x_1 + x_2)^2 (x_1 + x_2^2 +  x_2) ,       \;  \;     (x_1 - x_2)(x_1 + x_2)^2 (x_1^2 + x_2^2)  \}       \\
\label{eq:MWZ46}
 & \{ (x_1 - x_2)(x_1 + x_2) (x_1 + x_2^2 +  x_2) ,     \;\;          (x_1 - x_2)(x_1 + x_2) (x_1^2 + x_2^2)  \}   
\end{align}
\end{subequations}

\noindent
\textbf{System \eqref{eq:MWZ41}  for  \cite[Example 4.1]{MWZ:2012}}:
Our $\GIF$ algorithm \ref{alg:ProjInvBasis} with input tolerance $10^{-10}$ shows that the system is already in geometric
involutive form.  The corresponding Pommaret basis is given
in  \cite[Example 4.1]{MWZ:2012}.    The Pommaret basis looks different from the system, but is just
a linear combination of the system's polynomials to accomplish the Gr\"obner like requirement for its highest terms
under the term ordering prescribed in the problem.
The resulting coefficient matrix of this $\GIF$ form, is a full rank $m = 3$, $3 \times 10$ matrix which is input 
to the FDR algorithm.  Since it has rank $m = 3$, one facial reduction
yields a reduced $(10 - m) \times (10 - m) = 7 \times 7$ moment matrix.
Application of the FDR algorithm using the reduced moment matrix, yields convergence
 in 13 iterations and 0.09 secs, with a projected residual error of $10^{-14}$.
These statistics are shown in Table \ref{tab:FDR-GIF}.  The reduction in moment matrix size 
from $10 \times 10$ to a $7 \times 7$ matrix 
is recorded in the rightmost column of the Table by the fraction $\frac{10}{7}$.
Determination of this reduced moment matrix then yields the full $10 \times 10$ moment matrix of rank $r = 7$.
Since the dimension of the kernel for $\GIF$ form is $d = 7 = r$
Algorithm \ref{alg:GIF-MMtx} terminates with the input system as its output. 
It can be checked that the ideal generated by this system is real radical.
Our facial reduction algorithms in Section \ref{sect:facialimpl} provide checks for the existence of additional facial reductions.
They show that there are no additional facial reductions for this problem.

\noindent
\textbf{System \eqref{eq:MWZ44}  for  \cite[Example 4.4]{MWZ:2012}}:
This is very similar to the previous system \eqref{eq:MWZ41}. 
 As \cite{MWZ:2012} notes the coordinates for this example are not delta-regular, 
which they and we remedy by a linear change of  coordinates.  
We show that the original system is geometrically involutive, which is equivalent to the determination of a Pommaret basis by \cite{MWZ:2012}. 
Just as in the previous example, we form a $10 \times 10$ moment matrix from the $\GIF$ form, which is transformed by one facial reduction to a $7 \times 7$ matrix.  There are no additional facial reductions, and the full moment matrix and its rank $r$ are determined.
We find that dimension of the kernel for $\GIF$ form is $d = 7 = r$, so 
Algorithm \ref{alg:GIF-MMtx} terminates with the input system as its output. 
It can be verified the the output is a $\GIF$ form for the real radical of the ideal.

\noindent
\textbf{System \eqref{eq:MWZ42}  for  \cite[Example 4.2]{MWZ:2012}}:
This is quite similar to the systems \eqref{eq:MWZ42} and \eqref{eq:MWZ44}.
Our methods are similarly efficiently applied to this system.
Our $\GIF$ algorithm first applied one prolongation to the second system
\eqref{eq:MWZ42} to yield a degree $3$ system.
After projectiing from this degree $3$ system it shows that the resulting degree $2$ system is involutive and consists
of $3$ polynomials.
This degree $2$ system is geometrically equivalent to the Pommaret basis found by \cite{MWZ:2012}.  This system is simply the original
$2$ polynomials, together with their compatibility condition or S-polynomial 
 $x_2(x_1^2 - x_2 ) - x_1 (   x_1 x_2 - x_3 ) = x_1 x_3 - x_2^2 $.
Thus the input system $R$ is replaced with $ \ProjKer \ProKer  R$ with corresponding $3 \times 10$ coefficient matrix.
The resulting $10 \times 10$ moment matrix is facially reduced to a $7 \times 7$ moment matrix.
As in the previous examples, no new relations are detected in the kernel of the next moment matrix, $d=r=7$ and the algorithm
terminates.  It can be verified that the $\GIF$ form is a basis for the real radical ideal of the input system.

Unlike the systems \eqref{eq:MWZ41},\eqref{eq:MWZ42},\eqref{eq:MWZ44}, 
the remaining three systems
 \eqref{eq:MWZ43},\eqref{eq:MWZ45},\eqref{eq:MWZ46} of \cite{MWZ:2012} lead to new members in the kernel of their moment matrices.

\noindent
\textbf{System \eqref{eq:MWZ43}  for  \cite[Example 4.3]{MWZ:2012}}:
Our initial application of FDR showed slow convergence.
However a random linear change of coordinates applied to the input system $R$ dramatically improved the convergence.
Applying the $\GIF$ algorithm we found that $\Pro R$ is involutive and has a $8 \times 20$ coefficient matrix.
The dimension of its kernel is $d = 12$.
Facial reduction then reduces the $20 \times 20$ moment matrix to a $12 \times 12$ moment matrix
which has rank $r = 7 \not = d $ so the algorithm has not terminated.
The new member of the real radical arising in the moment matrix kernel can be alternatively derived by hand
by elimination of two of the systems polynomials:
$ x_1^2 + x_2^2 + x_3^2 -2  - ( x_1^2 + x_2^2 - x_3 ) = x_3^2 + x_3 - 2 =  (x_3 + 2)(x_3 - 1)$.
Then noting, as  explained in \cite{MWZ:2012}, that only the root $x_3 = 1$ leads to real solutions.
The $\GIF$  form of degree $2$ of the new system is computed.  Its coefficient matrix is $5 \times 10$
and has kernel of dimension $d =5$.
We note that even with the change of coordinates the FDR iteration of this second moment matrix did not initially converge
until we reduced the required projected residual error for production of the first moment matrix to $10^{-14}$.
The second moment matrix then was computed quickly and accurately as a $10 \times 10$ matrix which is reduced by one facial reduction to a $5 \times 5$ matrix. Since the rank of the moment matrix is $r = 5 = d$ our algorithm has terminated.
It can be checked that the output is equivalent to that found by \cite{MWZ:2012} and that the 
resulting $\GIF$ form is a basis for the real radical.

\noindent
\textbf{System \eqref{eq:MWZ45}  for  \cite[Example 4.5]{MWZ:2012}}:
Direct application of Algorithm \ref{alg:GIF-MMtx} to \eqref{eq:MWZ45} is relatively inefficient.
Instead  of this approach we consider an alternative subsystem approach which has the potential to be applied to larger systems.  Exploiting subsystem structure is a long established approach in system solving.

We apply Algorithm \ref{alg:GIF-MMtx} to the subsystem consisting of the first polynomial of  $P_1 = (x_1 - x_2)(x_1 + x_2)^2 (x_1 + x_2^2 +  x_2)$ of \eqref{eq:MWZ45}.
The $\GIF$ form of $P_1$ is just $P_1$, and its coefficient matrix 
is $1 \times 21 $ matrix with  a kernel of dimension $d = 20$.  
The corresponding moment matrix is $21 \times 21$, which is reduced
to a $ 20  \times 20  $ matrix after one facial reduction.   It has rank $r = 18  \not = d $.
So the algorithm has not terminated, and new members of the real radical are identified from the kernel of the moment matrix.
The new system is degree $ 5$ and has $3 $ polynomials. 
Algorithm $\GIF$ shows that the first projection of this system is involutive and is a single
fourth degree polynomial.
 Its coefficient matrix is $1  \times  15$  and its kernel has dimension $d = 14$.
The FDR algorithm produces a $ 15  \times 15 $ moment matrix which facially reduced to a
$14 \times 14$ moment matrix.  The rank of the moment matrix is $r = 14 = d$.
The algorithm terminates to coefficient errors
within $10^{-10}$ with output as a single polynomial which  is approximately:
\begin{equation}
\label{eq:RedP1}
 (x_1 - x_2)(x_1 + x_2) (x_1 + x_2^2 +  x_2)
\end{equation}
It can be checked that (\ref{eq:RedP1}) is a geometric involutive basis for the real radical for the ideal generated by $P_1$.

Similarly we apply Algorithm \ref{alg:GIF-MMtx} to the first polynomial of  \eqref{eq:MWZ45} which is given by $P_2 =  (x_1 - x_2)(x_1 + x_2)^2 (x_1^2 + x_2^2)$.
The algorithm now terminates with output as a single polynomial which is approximately:
\begin{equation}
\label{eqf}
 (x_1 - x_2)(x_1 + x_2) 
\end{equation}
This can be verified to be a geometric involutive  basis for the 
real radical for the ideal generated by $P_2$.

Then we consider the system
\begin{equation}
\label{eq:Ex45Redsys}
(x_1 - x_2)(x_1 + x_2) (x_1 + x_2^2 +  x_2), \; \; (x_1 - x_2)(x_1 + x_2) 
\end{equation}
The calculation for \eqref{eq:MWZ46}for Example 4.6 below yields a geometric involutive basis
which is approximately
\begin{equation}
\label{eq:Ex45sys}
(x_1^2 - x_2^2 )
\end{equation}
It can be independently checked that this is a $\GIF$ form for the real 
radical of the ideal of \eqref{eq:MWZ45}.

\noindent
\textbf{System \eqref{eq:MWZ46}  for  \cite[Example 4.6]{MWZ:2012}}:
This concerns the real solution of  $Q_1 = \eqref{eq:MWZ46} = \eqref{eq:MWZ46}$ subject to the constraints
$x_1 \geq 1$, $x_2 \geq 1$.  
Applying Algorithm \ref{alg:GIF-MMtx} to $Q_1$ yields a geometric involutive basis which 
is approximately $x_1^2 - x_2^2$.  This can be indepdently verified to be
a geometric basis for the real radical of $Q_1$.
The statistics of this reduction are given in the table
in the row labeled as Ex 4.6 $Q_1$.

To impose $x_1 \geq 1$, $x_2 \geq 1$ we substitute 
$ x_1 = x_3^2 + 1 ,  x_2 =  x_4^2 + 1$ and reduce the resulting polynomial $Q_2$
with Algorithm \ref{alg:GIF-MMtx}.  
We obtain $x_1 - x_2$ in agreement with  \cite[Example 4.6]{MWZ:2012}.
The statistics of this reduction are given in Table \ref{tab:FDR-GIF}
in the row labeled as Ex 4.6 $Q_2$.

\begin{table}
\tiny{
\setlength{\tabcolsep}{.56667em}
{\renewcommand{\arraystretch}{1.8}
\begin{tabular}{|lccccccc|}
  \hline
         &             & FDR   & FDR  & FDR           &GIF-FDR its & GIF &      Mom Mtx redn    \\
  Syst. & (n,d,p)  & \# its & secs & proj res err &    (\# FR )   & tol  &  factors $s(M)/s(\hat{M})$   \\
\hline
Ex4.1 &(3,2,3)    &  13    &  0.09 & $10^{-14}$&   1(1)         &$ 10^{-10}$ &  $\frac{10}{7}$   \\     
Ex4.2 &(3,2,2)    &  28     & 0.01 & $10^{-14}$&    1(1)         &$10^{-10}$ &  $\frac{10}{7}$     \\    
Ex4.3 &(3,2,2)    & 888, 238 & $2.3,0.6$  & $10^{-14}, 10^{-13}$&    2(2,1)        &$10^{-10}$ &   $\frac{20}{12},\frac{10}{5}$   \\     
Ex4.4 &(3,2,3)    & 346     &  0.53  & $10^{-14}$&    1(1)       & $10^{-10}$&   $\frac{10}{7}$    \\        
Ex4.5 $P_1$ &(2,5,1) &$22314,50 $  &37.6, 0.3& $10^{-12},10^{-14}$ &2 (2, 1) & $10^{-10}$& $\frac{21}{20},\frac{15}{14}$         \\    
Ex4.5 $P_2$ &(2,5,1) &$957,1 $  &4.4, 0.1 & $10^{-12},10^{-14}$ &2 (2, 1) & $10^{-10}$& $\frac{21}{20},\frac{6}{5}$         \\          
Ex4.6 $Q_1$ &(2,4,1)    & 170, 1   & 1.0, 0.09  &$10^{-12}, 10^{-14}$ &  2(2,1)        & $10^{-10}$& $\frac{21}{15},\frac{6}{5}$   \\   
Ex4.6 $Q_2$ &(1,4,1)    & 484, 1   & 1.4, 0.08  &$10^{-12}, 10^{-14}$ &  2(2,1)        & $10^{-10}$& $\frac{15}{14},\frac{6}{5}$   \\              
Cyl2d &(2,2,1)    &  10     & 0.19 &$10^{-15}$ &  1(1)          &$10^{-10}$ &  $\frac{6}{5}$   \\       
Cyl3d &(3,2,2)    &  33     & 0.77 &$10^{-14}$ &   1(1)         &$10^{-10}$ &  $ \frac{20}{12}$   \\       
Cyl4d &(4,2,3)    &  142   & 8.45 &$10^{-14}$ &  1(1)           &$10^{-10}$ & $\frac{70}{28}$    \\   
\hline  
\end{tabular}
}
}

\caption{\small{{\bf Statistics for the application of GIF and FDR to 
polynomial systems}:   $n = $ number of variables, $d = $ maximum polynomial degree, 
$p = $ the number of polynomials;
$s(M)$, $s(\hat{M})$ sizes of moment matrix $M$ and the faciallly
reduced matrix $\hat{M}$, resp.
Ex 4.1-4.6 are the $6$ examples in MWZ \cite{MWZ:2012}; Cyl2d-Cyl4d are 
the intersecting cylinder examples.
  }}
\label{tab:FDR-GIF}
\end{table}

\subsection{Intersecting higher dimensional cylinders}

Consider the systems of polynomials defining the intersection of $n - 1$ cylinders in $\mathbb{R}^n$
\begin{equation}
Cyl_{nd} := x_1^2 + x_2^2 - 1,   x_1^2 + x_3^2 - 1,  \cdots ,  x_1^2 + x_n^2 - 1 .
\end{equation}
Application of the $\GIF$ algorithm to the systems $Cyl_{nd}$ for $n = 2, 3, 4$ show that the systems become
geometrically involutive after $0, 2, 3$ prolongations respectively.
Table \ref{tab:FDR-GIF} shows the statistics for the subsequent application of Algorithm \ref{alg:GIF-MMtx} to these systems.
The algorithm converges quickly and accurately.
Indeed it can be independently determined that the it yields an geometric involutive basis for the real radical.

Further it can be determined that the cylinders form a complete intersection and the length of the prolongation
to make them involutive, can be determined from the symbol of the initial system \cite{MollerSauer:2000}.
The lower degree system, is geometrically formally integrable, and it would be interesting to develop
methods based on such lower degree systems, to determine, whether one can rule out new members
in the kernel of the moment matrix of the prolonged involutive system from such lower degree systems.

Finally we mention that recently certain so-called critical point methods have been developed 
for determining witness points \cite{WuReid13,Hauenstein12} on real components of real polynomial systems.
Indeed the method developed in \cite{WuReid13} is successful in finding a point on every component, if the ideal is both
real radical, and forms a regular sequence.  Consequently the systems above, the real radical is an important property for such solvers.
Such a regular sequence can be checked by dimension computation, 
we only need a formally integrable system which has lower degree than the involutive system, this leads to a smaller size of moment matrix.
Other interesting related results are given in \cite{MaZhi12}.

\subsection{Example of Matlab routine FDR}
\begin{example}
We first use the matrix from \eqref{eq:matB1}
\begin{equation}
\label{eq:matB1}
{B_1}^T = 
\begin{bmatrix}
     2   & 0 & 0 & 0&  -1    \\
  \end{bmatrix}.
\end{equation}
The moment matrix we get is the exactly the same as that in
\cite[Equation (37)]{ReidWangWu:14}:
\[
\begin{array}{rcl}
P
&=&
\begin{bmatrix}
 1.0000 &   -0.0000 &    1.4142 &   -0.0000 &    2.0000\cr
-0.0000 &    1.4142 &   -0.0000 &    2.0000 &   -0.0000\cr
 1.4142 &   -0.0000 &    2.0000 &   -0.0000 &    2.8284\cr
-0.0000 &    2.0000 &   -0.0000 &    2.8284 &   -0.0000\cr
 2.0000 &   -0.0000 &    2.8284 &   -0.0000 &    4.0000\cr
\end{bmatrix}
\end{array}
\]
The nullity/kernel matrix of $P$ is the same as in 
\cite[Equation (37)]{ReidWangWu:14} as well:
\[
\begin{bmatrix}
 2 &  0.026491 &  -0.23757 \cr 
 0 &  -0.81147 & -0.090484 \cr 
 0 &  -0.09366 &   0.83995 \cr 
 0 &   0.57379 &  0.063982 \cr 
-1 &  0.052982 &  -0.47515 \cr 
\end{bmatrix}
\]
though it is difficult to see from the last two columns.

To check whether the matrix $B_1$ in \eqref{eq:matB1} provides
the same nullity as the nullity of the matrix $P$,
one can look at the following short MATLAB code and see that it is so, 
i.e.,~the rank is correct and the spans do not change.
\begin{verbatim}
B1=[ B'
  sqrt2 0 -1 0 0
0 sqrt2 0 -1 0]
B1 = 
    2.0000         0         0         0   -1.0000
    1.4142         0   -1.0000         0         0
         0    1.4142         0   -1.0000         0

>> B1=B1'
B1 =

    2.0000    1.4142         0
         0         0    1.4142
         0   -1.0000         0
         0         0   -1.0000
   -1.0000         0         0
>> K=[null(P) B1]
K =
    0.8099    0.4053    0.1922    2.0000    1.4142         0
   -0.2574    0.1542    0.7593         0         0    1.4142
   -0.4913    0.6222   -0.2930         0   -1.0000         0
    0.1820   -0.1091   -0.5369         0         0   -1.0000
   -0.0575   -0.6426    0.1110   -1.0000         0         0
>> svd(K)
K>> svd(K)
    2.8284
    2.0000
    1.4142
    0.0000
    0.0000
\end{verbatim}

Following is the output during the MATLAB program. Note the quick and
accurate convergence; though we have to remember this is a tiny problem.
It took 118 iterations to get $15$ decimals accuracy. The moment matrix
$P$ has the correct rank.
\begin{verbatim}
Starting with new B value 
using [no*VV'] as initial starting point for P 
time for matrix repres.  0.0468003   
Starting while loop for Douglas-Rachford algorithm
iter       cos-vecs     norm-proj.-resid.   PSD-proj-per.iter.time 
 10           0.9938         0.04919         6.23e-05 
 20                1        0.005256        6.377e-05 
 30                1       0.0004443        6.188e-05 
 40                1       3.282e-05         6.23e-05 
 50                1       2.167e-06         0.000109 
 60                1       1.271e-07        6.467e-05 
 70                1        6.36e-09        6.551e-05 
 80                1       2.341e-10        6.251e-05 
 90                1       3.539e-12        6.349e-05 
 100                1       1.037e-12        6.439e-05 
 110                1       1.324e-13        6.572e-05 
 118                1       7.531e-15        6.404e-05 
time for iterations/while loop is  0.0780005   
max cosine value is  1  
checking feas error in DRalg.m using ***projected*** last iterate Rpsd 
error for norm(B'*P)  is  0 
\end{verbatim}

\end{example}

\section{Conclusion}
\label{sect:Conclusion}


SDP feasibility problems typically involve the intersection of the convex cone of semi-definite matrices with a linear manifold.  
Their importance in applications has led to the development of many specific algorithms.
However these feasibility problems are often marginally infeasible,
i.e.,~they do not satisfy strict feasibility as is the case for our polynomial applications. Such problems are
\emph{ill-posed} and \emph{ill-conditioned}. 


The main contribution of this paper is to introduce facial reduction, for the class of SDP problems arising
from analysis and solution of systems of real polynomial equations for real solutions.   
Facial reduction yields an equivalent problem for which there are strictly feasible points and which, in addition, are smaller.
Facial reduction also reduces the size of the moment matrices occurring in the application 
of SDP methods.  For example the determination of a $k \times k$ moment matrix for a problem
with $m$ linearly independent constraints is reduced to a $(k-m) \times (k-m)$ moment matrix by one facial reduction.  
We use facial reduction with our MATLAB implementation of Douglas-Rachford iteration (our FDR method).
In the case of only one constraint, say as in the case of univariate polynomials, one might expect that the improvement
in convergence due to that facial reduction would be minor.  However we present a class of geometric 
univariate polynomials of odd degree, where one such facial
reduction  combined with DR iteration, yields the real radical much more efficiently
than the standard interior point method Yalmip.
The high accuracy required by facial reduction and also the ill-conditioning commonly encountered in numerical
polynomial algebra \cite{Stetter:2004} motivated us to implement Douglas-Rachford iteration.


A fundamental open problem is to generalize the work of \cite{LasserreLaurentRostalski09,MR2830310} to positive dimensional ideals.
The algorithm of \cite{MWZ:2012,MA12} for a given input real polynomial system $P$, modulo the successful
application of SDP methods at each of its steps, computes a Pommaret basis $Q$:
\begin{equation}
\label{eq:RealRadConj}
\sqrt[\R]{ \left\langle P
\right\rangle_\R  }  \; \;   \supseteq  \; \;  \left\langle Q
\right\rangle_\R  \; \;  \supseteq  \; \;  \left\langle P
\right\rangle_\R 
\end{equation}
and would provided a solution to this open problem if it is proved that $ \left\langle Q
\right\rangle_\R = \sqrt[\R]{ \left\langle P
\right\rangle_\R  }$.
We believe that the work \cite{MWZ:2012,MA12} establishes an important 
feature -- involutivity -- that will necessarily be a 
a main condition of any theorem and algorithm characterizing the real radical.
Involutivity is a natural condition, since any solution of the above open problem using SDP, if it establishes radical
ideal membership, will necessarily need (at least implicitly) a real radical Gr\"obner basis.
Our algorithm, uses geometric involutivity, and similarly gives an intermediate ideal, which constitutes another 
variation on this family of conjectures.

In addition to implementing an algorithm to determine a first facial reduction.
We also implemented a test for the existence of additional facial reductions beyond the first
(e.g. in the cases of Examples 4.3 and 4.5 of \cite{MWZ:2012}).
By using the CVX package or Douglas-Rachford iteration to solve for the auxiliary problem, 
we can determine that if we need a second facial reduction by checking whether the 
optimal value of the auxiliary problem is close to $0$.
So far only moderate improvements in convergence have been obtained by our preliminary
implementation for construction of additional facial reductions.


Numerical polynomial algebra has been a rapidly expanding and popular area \cite{Stetter:2004}.
It's problems are typically very demanding, motivating the implementation of methods to improve accuracy.
For example Bertini, the homotopy package developed for numerical polynomial algebra, uses variable precision
arithmetic, with particularly demanding problems requiring thousands of digits of precision.
Consequently this is also a motivation to develop higher accuracy methods, such as the FDR method of this paper.
Manipulations with radical ideals would be a by-product from such work.

We provided a small set of examples, that illustrate some aspects of our algorithms.
In Maple all of our examples were executed with Maple's $Digits := 15$ and the input tolerance $:= 10^{-10}$  for the $\gif$ algorithm
whch intensively uses LAPack's SVD.
Accuracy in the projected residual error for our tests were between $10^{-14}$ and $10^{-12}$.
The normalized generators obtained for our experiments had coefficients differing less than
$10^{-10}$ from the exact coefficients.

Our implementation of auxiliary facial reductions, as still preliminary and needs improvement.
Even if the real radical is theoretically accessible, the conditioning of the polynomial system, as measured by
the sensitivity of changes in the solutions to changes in the coefficients, is a significant computational affect.
So a more detailed study of this aspect is worthwhile.

\cleardoublepage
\addcontentsline{toc}{section}{Index}
\label{ind:index}
\printindex

\bibliography{.master,.edm,.psd,.bjorBOOK,.GI}

\end{document}